\documentclass[12pt]{amsart}

\usepackage{amscd}
\usepackage{srcltx}
\usepackage{amsmath}
\usepackage{amsfonts}
\usepackage{verbatim}
\usepackage{amsthm}
\usepackage{amssymb, mathrsfs}
\usepackage{amscd}
\usepackage{latexsym}
\usepackage{array}
\usepackage{enumerate}
\usepackage{longtable}
\usepackage[all]{xypic}
\usepackage[latin1]{inputenc}

\ifx\pdfpagewidth\undefined 
\usepackage[T1]{fontenc}
\else 
\usepackage[OT1]{fontenc} 
\fi 
\usepackage{amsfonts,amssymb,latexsym,longtable,amsmath}

\def\obr{^{-1}}
\def\stm{\setminus}

\def\norm#1{\left\Vert#1\right\Vert}

\def\End{{\mathrm {End}}\,}
\def\Tr{{\mathrm {Tr}}\,}
\def\Re{{\mathrm {Re}}\,}
\def\dist{\hbox{dist\,}}
\def\Int{\hbox{Int\,}}

\newtheorem{thm}{Theorem}[section]
\newcommand{\bthm}{\begin{thm}} \newcommand{\ethm}{\end{thm}}
\newtheorem{prop}[thm]{Proposition}
\newcommand{\bprp}{\begin{prop}} \newcommand{\eprp}{\end{prop}}
\newtheorem{fact}[thm]{Fact}
\newcommand{\bfct}{\begin{fact}} \newcommand{\efct}{\end{fact}}
\newtheorem{prob}[thm]{Problem}
\newcommand{\bprb}{\begin{prob}} \newcommand{\eprb}{\end{prob}}
\newtheorem{question}[thm]{Question}
\newcommand{\bque}{\begin{question}} \newcommand{\eque}{\end{question}}
\newtheorem{lem}[thm]{Lemma}
\newcommand{\blem}{\begin{lem}} \newcommand{\elem}{\end{lem}}
\newtheorem{claim}[thm]{Claim}
\newcommand{\bclm}{\begin{claim}} \newcommand{\eclm}{\end{claim}}
\newtheorem{cor}[thm]{Corollary}
\newcommand{\bcor}{\begin{cor}} \newcommand{\ecor}{\end{cor}}
\newtheorem{conj}[thm]{Conjecture}
\newcommand{\bcnj}{\begin{conj}} \newcommand{\ecnj}{\end{conj}}
\theoremstyle{definition}
\newtheorem{defn}[thm]{Definition}
\newcommand{\bdfn}{\begin{defn}} \newcommand{\edfn}{\end{defn}}
\newtheorem{spec}[thm]{Specializing}
\newcommand{\bspc}{\begin{spec}} \newcommand{\espc}{\end{spec}}
\theoremstyle{remark}
\newtheorem{rem}[thm]{Remark}
\newcommand{\brem}{\begin{rem}} \newcommand{\erem}{\end{rem}}
\newtheorem{cnv}[thm]{Convention}
\newcommand{\bcnv}{\begin{cnv}} \newcommand{\ecnv}{\end{cnv}}
\newtheorem{exam}[thm]{Example}
\newcommand{\bexm}{\begin{exam}} \newcommand{\eexm}{\end{exam}}
\newcommand{\bpf}{\begin{proof}} \newcommand{\epf}{\end{proof}}

\newcommand{\C}{\mathbb C}

\newcommand{\N}{\mathbb N}

\renewcommand{\phi}{\varphi}
\renewcommand{\theta}{\vartheta}

\newcommand{\A}{{\rm a}}

\newcommand{\gep}{{\epsilon}}

\newcommand{\gd}{{\delta}}

\newcommand{\s}{{\sigma}}

\begin{document}
\title{Precompact groups and property (T)}

\author[M. Ferrer]{M. Ferrer}
\address{Universitat Jaume I, Instituto de Matem\'aticas de Castell\'on,
Campus de Riu Sec, 12071 Castell\'{o}n, Spain.}
\email{mferrer@mat.uji.es}

\author[S. Hern\'andez]{S. Hern\'andez}
\address{Universitat Jaume I, INIT and Departamento de Matem\'{a}ticas,
Campus de Riu Sec, 12071 Castell\'{o}n, Spain.}
\email{hernande@mat.uji.es}

\author[V. Uspenskij]{V. Uspenskij}
\address{Department of mathematics, 321 Morton Hall, Ohio University, Athens, Ohio 45701, USA}
\email{uspenski@ohio.edu}

\thanks{ The first and second listed
authors acknowledge partial support by the Spanish Ministry of
Science, grant MTM2008-04599/MTM; and by Fundaci\'o Caixa
Castell\'o (Bancaixa), grant P1.1B2008-26}

\begin{abstract}
For a topological group $G$ the dual object $\widehat G$ is defined as the set of equivalence classes
of irreducible unitary representations of $G$ equipped with the Fell topology. It is well known that
if $G$ is compact, $\widehat G$ is discrete.
In this paper, we investigate to what extent this remains true for precompact groups,
that is, dense subgroups of compact groups. We show that:
(a) if $G$ is a metrizable precompact group, then $\widehat G$ is discrete;
(b) if $G$ is a countable non-metrizable precompact group, then $\widehat G$ is not discrete; (c) every non-metrizable
compact group contains a dense subgroup $G$ for which $\widehat G$ is not discrete.
This extends to the non-Abelian case what was known for Abelian groups.
We also prove that if $G$ is  a countable Abelian precompact group, then $G$
does not have Kazhdan's property (T), although $\widehat G$ is discrete if $G$ is metrizable.
\end{abstract}



\date{20 November 2012}

\maketitle \setlength{\baselineskip}{24pt}
\setlength{\parindent}{1cm}


\section{Introduction}
\label{s:intro}

For a topological group $G$ let $\widehat G$ be the set of equivalence classes
of irreducible unitary representations of $G$. The set $\widehat G$ is equipped
with the so-called Fell topology  \cite{fel1}, which can be
defined on every set of equivalence classes of not necessarily irreducible representations
of a given topological group $G$ (see Section~\ref{s:prelim} for a definition).
We recall what it means in some familiar cases (see \cite{Robert}):
\begin{enumerate}
\item[(1)] if $G$ is Abelian, then $\widehat G$
is the standard Pontryagin-van Kampen dual group and the Fell topology
on $\widehat G$ is the usual compact-open topology;
\item[(2)] when $G$ is compact, the Fell topology
on $\widehat G$ is the discrete topology;
\item[(3)] when $\widehat G$ is neither Abelian nor compact,
$\widehat G$ usually is non-Hausdorff.
\end{enumerate}

In general, little is known about the properties
of the Fell topology. Its understanding heavily depends on Harmonic Analysis.
For example, if $G$ is a second countable locally compact group, then
the Fell topology on $\widehat G$ satisfies
the $T_1$ separation axiom
if and only if $G$ is type I (see \cite{bhv:book}).

It is possible
to define a different topology on $\widehat G$ that, when $G$ is compact, is \emph{grosso modo}
the natural quotient of the set of all irreducible representations equipped with the compact-open topology.
However, this topology is less useful than the Fell topology because for every integer $n$ it makes the set of all
$n$-dimensional
representations closed, while it is often desirable to regard certain lower-dimensional representations
as limits of higher-dimensional ones (see \cite{bhv:book,dixmier,harp_vale,Marg,Robert}).

A topological group $G$ is {\em precompact}
if it is isomorphic (as a topological group) to a subgroup of a compact
group $H$ (we may assume that $G$ is dense in $H$).
If $G$ is compact, 
then $\widehat G$ is discrete. If $G$ is a dense
subgroup of $H$, the natural mapping $\widehat H\to \widehat G$ is a bijection but
in general need not be a homeomorphism.
Following Comfort, Raczkowski
and Trigos-Arrieta \cite{comractri:04}, we say that  $G$ {\em determines} $H$ if $\widehat G$ is discrete
(equivalently, if the natural bijection $\widehat H\to \widehat G$ is a homeomorhism).
 A compact group $H$ is {\em determined}
if every dense subgroup of $H$ determines $H$.
In the present paper we investigate the following
question: if $H$ is a compact group and $G$ is a dense subgroup of $H$, under what conditions
does $G$ determine $H$?
Equivalently, for what
precompact groups $G$ is $\widehat G$ discrete?

In the Abelian case, this question has been settled in the work of several authors.
Aussenhofer \cite{aus} and, independently, Chasco
\cite{chasco} showed that every metrizable Abelian compact
group  $H$ is  determined.
Comfort, Raczkowski and Trigos-Arrieta
\cite{comractri:04} noted that the Aussenhofer-Chasco theorem fails
for non-metrizable  Abelian compact groups $H$.
More precisely, they proved that every
non-metrizable compact Abelian group $H$ of weight $\geq 2^\omega$
contains a dense subgroup that does not determine $H$. Hence,
under the assumption of the continuum hypothesis, every determined
compact Abelian group $H$ is metrizable. Subsequently, it was shown in
\cite{hermactri} that the result also holds without assuming the
continuum hypothesis (see also \cite{dik_sha:jmaa_det}).

Our goal is to extend the results quoted above to compact groups that are not necessarily Abelian.
We now formulate our main results.

{\bf Theorem~\ref{t:1}}. {\em If $G$ is a precompact metrizable group, then $\widehat G$ is discrete.
Equivalently, every metrizable compact group is determined.}

This extends the aforementioned results from \cite{aus} and \cite{chasco} to non-Abelian
compact metrizable groups.

A certain extension of the Aussenhofer-Chasco theorem to non-Abelian groups is due to
Luk\'acs in \cite{lukacs:determined}, where he considers the compact open topology on certain groups of
continuous mappings. Theorem~\ref{t:1} does not follow from Luk\'acs' results because
in general the Fell topology on the dual object is not reduced to the compact-open topology.

In this paper ``countable" means ``countable infinite". Let $1_G$ be the class of the trivial representation.

{\bf Theorem~\ref{t:2}}. {\em If $G$ is a countable precompact non-metrizable group, then $1_G$ is not an isolated point in $\widehat G$.}

Theorem~\ref{t:2} also shows that {\em  if $G$ is a countable dense subgroup of a non-metrizable compact group $H$,
then $G$ does not determine $H$.}

{\bf Theorem~\ref{t:3}}. {\em If $H$ is a non-metrizable compact group, then $H$ has a dense subgroup $G$
such that $\widehat G$ is not discrete.}

Together with  Theorem~\ref{t:1} this shows that {\em a compact group is determined if and only if
it is metrizable.} This extends to non-Abelian compact groups the results given in
\cite{comractri:04,hermactri,dik_sha:jmaa_det} for Abelian compact groups.

The group $G$ has {\em property (T)} if $1_G$ is isolated in $\mathcal R \cup \{1_G\}$ for every set $\mathcal R$
of equivalence classes of unitary representations of $G$ without non-zero invariant vectors.
This definition is equivalent to the definition of property (T)
in terms of Kazhdan pairs which we recall below in Section~\ref{s:prelim} (see Proposition 1.2.3 in  \cite{bhv:book}).
Compact groups have property (T), and we are interested in property (T) for precompact groups.

It follows from Theorem~\ref{t:2} that every countable precompact group with property (T) is metrizable.
Furthermore, our proof also yields the following result of Wang
\cite{wang:1975}: if a discrete group $G$ has property (T), then its Bohr compactification $bG$ is metrizable.
(see Corollary \ref{co:wang} below).

{\bf Theorem~\ref{t:4}}. {\em If $G$ is a countable Abelian precompact group, then $G$ does not have property (T).}

We do not know whether there exists a non-compact precompact Abelian group with property (T) (Question~\ref{q:1}).
According to Theorem~\ref{t:4}, such a group must be uncountable.

\section{Preliminaries: Fell topologies and property (T)}
\label{s:prelim}

All topological groups are assumed to be Hausdorff.
For a (complex) Hilbert space $\mathcal{H}$ the unitary group $U(\mathcal{H})$
of all linear isometries of $\mathcal{H}$ is equipped with the strong operator topology
(this is the topology of pointwise convergence). With this topology, $U(\mathcal{H})$  is a topological group.
If $\mathcal{H}=\C^n$, we identify $U(\mathcal{H})$ with the {\em unitary group $\mathbb{U}(n)$
of order $n$}, that is, the compact Lie group of all complex $n\times n$ matrices $M$
for which $M^{-1}=M^{*}$.

A {\em unitary representation} $\rho$ of the to\-po\-lo\-gi\-cal group $G$
is a continuous homomorphism $G\to U(\mathcal{H})$, where $\mathcal{H}$ is a complex Hilbert space.
A closed linear subspace $E\subseteq \mathcal H$ is an \emph{invariant} subspace
for $\mathcal S\subseteq U(\mathcal{H})$ if $ME\subseteq E$ for all
$M\in \mathcal S$.
If there is a closed subspace $E$ with $\{0\}\subsetneq E\subsetneq \mathcal H$
which is invariant for $\mathcal S$, then $\mathcal S$ is called
\emph{reducible}; otherwise $\mathcal S$ is \emph{irreducible}.
An \emph{irreducible representation} of $G$ is a 
unitary representation $\rho$ such that 
$\rho(G)$ is irreducible.

Two unitary representations $\rho:G\to U(\mathcal{H}_1)$ and $\psi: G\to U(\mathcal{H}_2)$
are {\it equivalent} 
if there exists a Hilbert space isomorphism $M:\mathcal {H}_1\to \mathcal H_2$
such that $\rho(x)=M^{-1}\psi(x)M$ for all
$x\in G$. The {\em dual object} of a topological group $G$ is the set $\widehat G$ of
equivalence classes of irreducible unitary representations of $G$.

If $G$ is a precompact group, the Peter-Weyl Theorem (see \cite{hof_mor:compact_groups}) implies that
all irreducible unitary representation of $G$ are finite-dimensional and determine an embedding of $G$ into the product
of unitary groups $\mathbb{U}(n)$.

If $\rho:G\to U(\mathcal{H})$ is a unitary representation, a complex-valued function $f$ on $G$ is called
a {\em function of positive type} (or {\em positive-definite function}) {\em associated with} $\rho$ if there exists a vector
$v\in \mathcal{H}$ such that $f(g)=(\rho(g)v, v)$
(here $(\cdot ,\cdot)$ denotes the inner product in $\mathcal{H}$).  We denote by
$P_\rho'$ be the set of all functions of positive type associated with $\rho$.
Let $P_\rho$ be the convex cone generated by $P_\rho'$, that is, the set of sums of elements of $P_\rho'$.
Observe that if $\rho_1$ and $\rho_2$
are equivalent representations, then $P'_{\rho_1}=P'_{\rho_2}$ and $P_{\rho_1}=P_{\rho_2}$.

Let $G$ be a topological group, $\mathcal R$ a set of equivalence classes of unitary representations of $G$.
The \emph{Fell topology} on $\mathcal R$ is defined as follows: a typical neighborhood of $[\rho]\in \mathcal R$ has the form
$$
W(f_1, \cdots, f_n, C, \gep)=\{[\s]\in\mathcal R :
\, \exists g_1, \cdots, g_n\in P_\s
\  \forall x\in C  \
|f_i(x)-g_i(x)| <\gep\},
$$
where $f_1, \cdots, f_n\in P_\rho$ (or $P_\rho'$), $C$ is a compact subspace of $G$, and $\gep>0$.
In particular, the Fell topology is defined on the dual object $\widehat G$. Even though we will not use
the Jacobson topology in this paper, it is perhaps worth mentioning that, if $G$ is locally compact,
the Fell topology on $\widehat G$ can be derived from the Jacobson topology on the primitive ideal space of $C^*(G)$,
the $C^*$-algebra of $G$ \cite[section 18]{dixmier}, \cite[Remark F.4.5]{bhv:book}.

As mentioned in the Introduction, the group $G$ has {\em property (T)} if the trivial representation
$1_G$ is isolated in $\mathcal R \cup \{1_G\}$ for every set $\mathcal R$
of equivalence classes of unitary representations of $G$ without non-zero invariant vectors.

Let $\pi$ be a unitary representation of a topological group $G$ on a Hilbert space $\mathcal H$,
$F\subseteq G$, and $\gep>0$.
A unit vector $v\in \mathcal H$
is called $(F,\gep)$-\emph{invariant} if
$\norm{\pi(g)v-v}<\gep$ for every $g\in F$.
A topological group $G$ has property (T) if and only if there is compact subset $Q$ of $G$ and $\gep>0$
such that for every unitary representation $\rho$ with a $(Q,\gep)$-invariant vector has a non-zero
invariant vector \cite[Proposition 1.2.3]{bhv:book}. The pair $(Q,\gep)$ is often referred to as a {\em Kazhdan pair}.

We refer to Fell's papers \cite{fel1,fel2}, the classical text by Dixmier \cite{dixmier} and the recent monographs by
de la Harpe and Valette \cite{harp_vale}, and Bekka, de la Harpe and Valette
\cite{bhv:book} for basic definitions and results concerning Fell topologies and property (T).

It is well known that every compact group $H$ has a unique normalized Borel regular invariant measure, the
Haar measure. We let $L^2(H)$ denote the Hilbert space that is constructed with the aid of this measure.
If $H$ is metrizable, then it follows from the Peter--Weyl Theorem \cite{hof_mor:compact_groups,Robert} that, up to equivalence,
there are countably many irreducible unitary representations of $H$. We henceforth enumerate them as $\rho_i$, $i\in \mathbb{N}$,
let $\chi_i= \Tr(\rho_i)$ be their characters and $d_i$ their dimensions. Put $P_i'=P_{\rho_i}'$
be the corresponding set of functions of positive type, let $P_i$ be the convex cone generated by $P_i'$,
and let $Q_i\subseteq C(H)$ be the linear space generated by $P_i$. Because $H$ is compact, the spaces
$Q_i$ are finite-dimensional ($\dim Q_i=d_i^2$) and are pairwise orthogonal in the Hilbert space $L^2(H)$.
Let $N_i=\{f\in P_i:f(e)=1\}$ be the space of normalized functions in $P_i$.
This is a compact subset of $Q_i$, 
(see Lemma~\ref{l:referee} below). Set $h_i=\chi_i/d_i\in N_i$ be the normalized character.
A basic and comprehensive reference for all notions and results mentioned here is the monograph by
Alain Robert \cite[Part~I]{Robert}.

\section{Precompact  groups }

In this section, we prove some general results about precompact groups that are essential
for proving one of our main results (Theorem~\ref{t:1}) in the next section.
We assume that $G$ as a dense subgroup of a compact group $H$,
and use the notation $\rho_i$, $P_i$, etc., introduced in the end of the previous section.
We assumed there that $H$ is metrizable, and the index $i$ runs over the set $\mathbb{N}$ of positive integers;
the results of this section remain true if we allow $H$ to be non-metrizable and accordingly allow $i$ to run
over an uncountable index set I.

Integrals over $H$ will be taken with respect to the normalized Haar measure on $H$.
If $X$ is compact, we consider the space $C(X)$ of continuous functions as a metric space, with the metric
defined by the sup-norm.

\blem\label{l:referee} For every $i\in I$, the convex subset $N_i$ of the finite-dimensional vector space $Q_i$
is compact.
\elem

Of course, ``compact" refers to the natural topology of the finite-dimensional vector space $Q_i$.
For greater clarity, the space $Q_i$ is equipped with the topology of $L^2(H)$, which coincides
on $Q_i$ with the topology inherited from $C(H)$, because $Q_i$ is finite dimensional \cite[Theorem~1.21]{rudin}.

\begin{proof}
Consider the set $N_i'=\{f\in P_i': f(e)=1\}$. We claim that $N_i$ is the convex hull of $N_i'$.
Indeed, given $f\in N_i$, write $f$ as a finite sum $f=\sum f_k$, where $f_k\in P_i'$. We may assume
that $f_k(e)\ne0$ for each $k$, since otherwise $f_k$ is identically zero and can be omitted.
Put $g_k=f_k/f_k(e)$.
Then $\sum f_k(e)=1$, and $f$ is the convex combination $f=\sum f_k(e)g_k$ of the functions $g_k\in N_i'$.

Denote by $S$ the unit sphere (= the set of all vectors of length one) of the space of the representation $\rho_i$.
Define a continuous onto mapping $S\to N_i'$ by sending every $v\in S$ to the function
$f_v\in N_i'$ defined by $f_v(g)=(\rho_i(g)v, v)$. Since $S$ is compact, so is $N_i'$.
It is well known that the convex hull of a compact set in a finite-dimensional vector space
is compact \cite[Theorem 3.20(d)]{rudin}. Thus $N_i$ is compact.
\end{proof}

Recall that a function on $G$ is {\em central} if it is constant on conjugacy classes
(cf. \cite{Robert}). The next lemma is a straightforward consequence of the properties
of characters in a compact group (see \cite[7.6]{Robert}). We include its proof here for the reader's convenience.

\blem \label{l:central}
The only central functions in $Q_i$ are the functions $c\chi_i$, $c\in \C$.
\elem
\begin{proof} Any central function $g\in Q_i$ is the sum in $L^2(H)$ of the series $\sum c_j\chi_j$, where $c_j=\int_H g\bar{\chi_j}$
\cite[7.6]{Robert}. Since $Q_i\perp Q_j$, we have $c_j=0$ for all $j\ne i$.
\end{proof}

\blem \label{l:dist}
Let $X$ be compact, $D$ a dense subset of $X$ and  $N$ a compact subset of $C(X)$.
If $g\in C(X)$ is at the distance $>\gep$
from $N$, there exists a finite subset $F\subseteq D$ such that the distance from
$g|_F$ to $N|_F$ in $C(F)$ is $>\gep$.
\elem

\begin{proof}
For every $f\in N$ there exists a point $x=x(f)\in D$ such that $|f(x)-g(x)|>\gep$. Pick a neighborhood
$O_f$ of $f$ such that $|h(x)-g(x)|>\gep$ for every $h\in O_f$. Since $N$ is compact, it is covered by a finite
collection of such neighborhoods, say $O_{f_1}, \dots, O_{f_s}$. The set $F=\{x(f_1), \dots, x(f_s)\}$
is as required.
\end{proof}

\blem\label{l:separ}
The space $\widehat G$, equipped with the Fell topology, is $T_1$.
\elem

\begin{proof}
Let $\rho_i$ and $\rho_j$ be non-equivalent irreducible unitary representations of $G$.
We construct a neighborhood of $[\rho_i]$ in $\widehat G$ which does not contain $[\rho_j]$.
Consider the set $K_j=\{f\in P_j: f(e)\le 2\}$. This set is compact in the natural topology
of the finite-dimensional space $Q_j$ that contains $P_j$ and $K_j$. Indeed, $K_j$ is the image
of $[0,2]\times N_j$ under the mapping $(t,f)\mapsto tf$ ($0\le t\le2$, $f\in N_j$), and $N_j$ is compact
(Lemma~\ref{l:referee}). The normalized character $h_i=\chi_i/d_i$ of $\rho_i$ is orthogonal to $Q_j$,
hence $\dist(h_i, K_j)\ge\dist(h_i, Q_j)>0$. It follows from Lemma~\ref{l:dist} that there exist a finite
set $F\subseteq G$ and $\gep>0$ such that $\dist(h_i|_F, K_j|_F) >\gep$. We may assume that $e\in F$ and
$\gep\le1$. We claim that the neighborhood $W(h_i, F, \gep)$ of $[\rho_i]$ does not contain $[\rho_j]$.
Equivalently, if $f\in P_j$, then $\dist(h_i|_F, f|_F)\ge \gep$. Indeed, we have just seen that this is true
if $f\in K_j$. If $f\in P_j\stm K_j$, we have $f(e)>2$ and hence
$\dist(h_i|_F, f|_F)\ge |f(e)-h_i(e)|=f(e)-1>1\ge\gep$.
\end{proof}

The next lemma is a straightforward consequence of the classical
Riemann-Lebesgue lemma for compact groups (see \cite[28.40]{hr:ii}). Again, we include its proof here for the reader's
convenience. If $\{x_i: i\in I\}$ is a family of numbers, the notation $x_i\to 0$ means that for every $\gep>0$
the set $\{i\in I: |x_i|>\gep\}$ is finite.

\blem \label{l:int}
Let $V$ be a measurable subset of $H$. 
Then $\int_V \chi_i\to 0$ as $i\in I$.
\elem
\begin{proof}
The integrals $\int_V \chi_i$ are the 
scalar products of the characteristic function of $V$ with the terms of the orthonormal
sequence $(\bar\chi_i)$ in $L^2(H)$.
\end{proof}

\blem\label{l:new}
Let $V\subseteq H$ be a compact neighborhood of the identity $e$ that is invariant under inner automorphisms of $H$.
Let $f\in C(H)$ be a continuous central function.
Then the distance from $f|_V$ to $N_i|_V$  in the space $C(V)$ is attained at $h_i|_V$,
where $h_i=\chi_i/d_i\in N_i$.
\elem
\begin{proof}
Let $h\in N_i$ be such that $h|_V$ is as close to $f|_V$ as possible. Here we use the compactness of $N_i$,
see Lemma~\ref{l:referee}. If $g$ is a function on $H$ and $y\in H$, let $g_y$ be the function on $H$ defined by
$g_y(x)=g(yxy\obr)$. The formula $(\rho_i(x)v,v)=(\rho_i(yxy\obr)\rho_i(y)v, \rho_i(y)v)$, where $v$ is any vector
in the space of the representation $\rho_i$, shows that $P_i'$ is invariant under the mapping $h\mapsto h_y$,
hence $N_i$ is invariant as well. Let $h'$ be the function on $H$ defined by $h'(x)=\int_{y\in H}h(yxy\obr)$;
equivalently, $h'=\int_{y\in H}h_y$, where we use vector-valued integrals, as in \cite[Theorem 3.27]{rudin}.
According to the cited theorem, $h'\in N_i$. For every $y\in H$
the distance from $h|_V$ to $f|_V$ is the same as the distance from $h_y|_V$ to $f_y|_V$. Since $f$ is central,
we have $f_y=f$. Thus all $h_y|V$ lie in the compact convex set of points of $N_i|_V$ that are as close as possible
to $f|_V$. It follows that $h'|_V=\int_{y\in H}h_y$ also lies in this set. Therefore the distance from $f|_V$ to $N_i|_V$
is attained at $h'|_V$.
It is easy to verify that $h'_y=h'$ for every $y\in H$, which means that $h'$ is a central function.
By Lemma~\ref{l:central}, the only central function in $N_i$ is $h_i=\chi_i/d_i$. Thus $h'=h_i$.

%
\end{proof}

\section{Precompact  metrizable groups}

The aim of this section is to prove the following result.
\bthm \label{t:1}
If $G$ is a precompact metrizable group, then $\widehat G$ is discrete.
\ethm

The idea of the proof of Theorem \ref{t:1} is as follows. We want to show that every point $[\rho]\in\widehat G$
is isolated. Since  $\widehat G$ is $T_1$ (Lemma~\ref{l:separ}), it suffices to find
a neighborhood $W$ of $[\rho]$ which for some integer $i_0$ does not contain any $[\rho_i]$
with $i\ge i_0$. Our neighborhood will be of the form $W=W(h,F,\gep)$, where $h$ is the normalized
character of $[\rho]$ and $F=\{e\}\cup\bigcup_{i\ge i_0} F_i$ is a compact subset of $G$, where $(F_i)$ is
as sequence of finite sets which converges to $e$ in the sense that every neighborhood of $e$
contains every $F_i$ with a sufficiently large index. The finite set $F_i$ ``takes care" of $\rho_i$,
in the sense that
it ensures that the neighborhood $W$ does not contain $[\rho_i]$. We derive the existence
of $F_i$ from the orthogonality of characters. If $V$ is a neighborhood of $e$ on which $h$ is
close to $1$, the orthogonality relations imply that $\int_V \chi_i \to0$ as $i\to\infty$ (Lemma~\ref{l:int}),
which forces $\chi_i$ to be close to $0$ somewhere on $V$ for $i\ge i_0$. This implies that $h$ and $h_i$
are not close to each other on $V$. With a little more work we show that $h$ is not close to any element
of $P_i$, and that this is witnessed by a certain finite subset $F_i$ of $V$.

We now are ready to prove Theorem~\ref{t:1}.

\bpf
Recall that we view $G$ as a dense subgroup of a compact metrizable group $H$.
Pick a bi-invariant (i.e., invariant under left and right translations, hence also under inner automorphisms) metric $b$ on $H$.
We denote by $V_\gep$ the closed $\gep$-ball (with respect to $b$) centered at the neutral element $e$.
We use the notation introduced above.

We can identify the sets $\widehat G$ and $\widehat H$ (not taking the topology of these sets into account).
Let $\rho$ be an irreducible unitary representation of $H$. We show that $[\rho]$ is an isolated point in $\widehat G$.

Let $\chi=\chi_\rho$ be the character of $\rho$, and $d=d_\rho$ its dimension. Consider the normalized character
$h=\chi/d$.
Pick $\gep>0$  such that $\Re h(x)>2/3$ for every $x\in V_\gep$.
Put $V=V_\gep$. Since $\int_V \chi_i\to 0$ (Lemma~\ref{l:int}), there is an index $i_0$ such that for all
$i\geq i_0$ there exists a point $x_i\in V$
such that $\Re \chi_i(x_i)\le 1/3$. Choose such an $x_i$ as close to $e$ as possible. In other words,
for $i\ge i_0$ the point $x_i$ is a point in the compact set $\{x\in V: \Re \chi_i(x_i)\le 1/3\}$ closest to $e$.


 Let  $\gep_i=b(x_i, e)$. We claim that $x_i\to e$ in the metric $b$ (equivalently, that $\gep_i\to 0$). Indeed, otherwise there exist a $\gd>0$
and an infinite set $A$ of integers $\ge i_0$ such that $\gep_i>\gd$ for every $i\in A$. Then $\Re \chi_i(x)>1/3$ for every
$x\in V_\gd$ and every $i\in A$, since $x$ is closer to $e$ than $x_i$. This contradicts the fact that
$\int_{V_\gd} \chi_i\to 0$ (Lemma~\ref{l:int}).

As before, let $h_i=\chi_i/d_i$.
Pick $\gep_i'>\gep_i$ so that $\gep_i'\to 0$.
Put $V_i=V_{\gep_i}$, $V_i'=V_{\gep_i'}$. Note that $\Re h(x_i)>2/3$ and
$\Re h_i(x_i)=\Re\chi_i(x_i)/d_i\le1/3$, so $|h(x_i)-h_i(x_i)|>1/3$ and hence
$\hbox{dist\,}(h|_{V_i}, h_i|_{V_i})>1/3$. It follows from Lemma~\ref{l:new}
that $\hbox{dist\,}(h|_{V_i}, N_i|_{V_i})>1/3$.
We claim that there exists a finite set $F_i\subseteq V_i'\cap G$ such that
$\hbox{dist\,}(h|_{F_i}, N_i|_{F_i})>1/3$.
This follows from Lemma~\ref{l:dist}, in which we make the following identifications:
using the notation of  Lemma~\ref{l:dist}, we set $D=G\cap\Int V_i'$, where Int denotes the interior
with respect to $H$; for $X$ we take the closure of $D$ in $H$, which is the same as the closure
of $\Int V_i'$, since $G$ is dense in $H$;
and the role of the compact set $N\subseteq C(X)$ is now played by $N_i|_X$, where $X$ is as above.
The compactness of $N_i$ was established in Lemma~\ref{l:referee}.
Since $V_i\subseteq \{x\in H:b(x,e)<\gep_i'\}\subseteq \Int V_i'\subseteq X$, we have
$\hbox{dist\,}(h|_X, N_i|_X)\ge\hbox{dist\,}(h|_{V_i}, N_i|_{V_i})>1/3$,
so the conditions of Lemma~\ref{l:dist} are satisfied with $\gep=1/3$.

Set $F=\{e\}\cup\bigcup\limits_{i\geq i_0} F_i$. Since $\gep_i'\to 0$,  $F$ is a compact subset of $G$.
Observe that $\hbox{dist\,}(h|_F, N_i|_F)\ge \hbox{dist\,}(h|_{F_i}, N_i|_{F_i})>1/3$ for every $i\ge i_0$.
We claim that the neighborhood $W(h|_G, F, 1/6)$ of $[\rho]$ does not contain any $[\rho_i]$ with $i\ge i_0$.
Indeed, let $i\ge i_0$ and $f\in P_i$. We show that $\hbox{dist\,}(h|_F, f|_F)\ge 1/6$. Assume that
 $\hbox{dist\,}(h|_F, f|_F)<1/6$. Put $c=f(e)$ and $g=f/c\in N_i$.
Since $e\in F$, we have  $|1-c|=|h(e)-f(e)|<1/6$.
Since $g$ is a function of positive type,
we have $|g(x)|\le g(e)=1$ for every $x\in H$  \cite[Proposition C.4.2(ii)]{bhv:book}.
It follows that
$|f(x)-g(x)|=|cg(x)-g(x)|=|c-1|\cdot|g(x)|<1/6$ for every $x\in H$, hence $\hbox{dist\,}(h|_F, g|_F)\le
\hbox{dist\,}(h|_F, f|_F)+\hbox{dist\,}(f|_F, g|_F)<1/6+1/6=1/3$. Since $g\in N_i$, this contradicts
the inequality $\hbox{dist\,}(h|_F, N_i|_F)>1/3$.

We have proved that the neighborhood $W(h|_G, F, 1/6)$ of $[\rho]$ contains only finitely many elements
of $\widehat G$. Since
$\widehat G$ is a T$_1$-space (Lemma~\ref{l:separ}), it follows that $\widehat G$ is discrete.
\epf

The compact set $F=F_\rho\subseteq G$ that we constructed in the proof above is the set of points of a sequence converging to $e$.
It depends on the point $[\rho]\in \widehat G$ the neighborhood of which we are constructing. However, with some little more
work, we can obtain a sequence which does not depend on the point we select in $\widehat G$.

\bcor
There exists a single compact subset $Q\subseteq G$ such that for every $[\rho_i]\in \widehat G$
the neighborhood $W(h_i|_G, Q, 1/6)$ of $[\rho_i]$ in $\widehat G$ is finite.
\ecor

\bpf
Indeed, in our definition of $F_\rho$ as the union $\{e\}\cup \bigcup_{i\ge i_0} F_i$ we clearly could have replaced
the index $i_0$ by any larger integer, without losing the property that $W(h|_G, F_\rho, 1/6)$ is finite.
Therefore we may assume that the diameter of $F_\rho$ is as small as we wish. Construct the sets $F_{\rho_1}, F_{\rho_2},\dots$
so that their diameters tend to~$0$, and let $Q$ be their union. If $U$ is any neighborhood of $e$, each set $F_{\rho_i}\stm U$ is finite,
and only finitely many of them are non-empty. It follows that
$Q$ is the set of points of a sequence converging to $e$ and hence compact.
Clearly $Q$ has the required property: for every $[\rho_i]\in \widehat G$
the neighborhood $W(h_i|_G, Q, 1/6)$ of $[\rho_i]$ in $\widehat G$ is finite,
being a subset of the finite set $W(h_i|_G, F_{\rho_i}, 1/6)$.
\epf

\brem
The referee suggested the following questions:
\begin{enumerate}
\item[(i)] For a precompact group $G$, does the Fell topology on $\widehat G$
contain the co-countable topology? That is, is every countable subset closed?
\item[(ii)] What if $G$ is pseudocompact?
\item[(iii)] What if $G$ contains no infinite compact subsets?
\item[(iv)] What if $G$ contains many compact $G_\delta$-subgroups?
\end{enumerate}
\erem
\medskip

Regarding the first three questions, let $G$ be an Abelian pseudocompact group without infinite compact
subsets (see \cite{AChDT,galmac:2011} and the references there).
The dual $\widehat G$ coincides with
the ordinary Pontryagin dual, which is a topological group. Furthermore, since $G$ has no infinite
compact subsets, the compact open topology on $\widehat G$ coincides with the pointwise convergence
topology. It follows that $\widehat G$ is a subgroup of a power of the circle group and hence is precompact.
Let $H$ be the completion of $\widehat G$. Then $H$ is a compact group. If $F$ is a countable subset of $\widehat G$, then there is a point
$p\in H$ that is an accumulation point of $F$. Therefore the identity of $\widehat G$ is an
accumulation point of the countable set $A=\{xy^{-1} : x,y\in F,\ x\not=y\}$. Since
the neutral element does not belong to $A$, this yields a countable subset of $\widehat G$ which
is not closed. Therefore the answer to the referee's first three questions is ``not always".

Regarding the fourth question, Theorem \ref{t:2} shows that the answer is also ``not always''.
Indeed, in a countable group every subgroup is a $G_\delta$-set. Furthermore, we can take the product
of a compact metrizable group and a countable precompact group of uncountable character in order to obtain
a precompact group with many compact $G_\delta$-subgroups such that its dual is not metrizable.
On the other hand, Theorem \ref{t:1},  asserting that if $G$ is precompact metrizable group
then $\widehat G$ is discrete, remains true for almost metrizable groups (that is, groups $G$ containing
a compact subgroup $K$ such that $G/K$ is a metrizable space). The proof extends the
ideas used here but full details will appear elsewhere.

\section{Countable precompact groups}

In this section we prove the following:

\bthm
\label{t:2}
If $G$ is a countable precompact non-metrizable group, then $1_G$ is not an isolated point in $\widehat G$.
\ethm

\bthm
\label{t:3}
If $H$ is a non-metrizable compact group, then $H$ has a dense subgroup $G$ such that $\widehat G$ is not discrete.
\ethm

For a topological group $G$ let $\widehat G_n\subseteq \widehat G$ be the set of classes of $n$-dimensional irreducible
unitary representations. 
We denote by $w(X)$ the weight
(i.e., the least cardinality of a base) of a topological space $X$. We allow the weight to be finite:
if $X$ a finite discrete space, $w(X)$ equals the cardinality of $X$.
We first establish two lemmas that will be used in the proof of Proposition~\ref{p:weight}.

\blem
\label{l:action}
If a compact group $H$ acts on a metric space $(M,d)$ by isometries, there is a natural metric on the quotient space $M/H$
compatible with the quotient topology, and $w(M/H)\le w(M)$.
\elem

\bpf
Let $Hx$ and $Hy$ be two $H$-orbits, where $x,y\in M$. Set the distance between $Hx$ and $Hy$ (considered as points in $M/H$) to be
the same as the distance between compact subsets $Hx$ and $Hy$ of the metric space $M$, which is $\min\{d(x',y'): x'\in Hx,\,y'\in Hy\}$.
This number is also equal to $\min\{d(x',y): x'\in Hx\}$ and to $\{d(x,y'): y'\in Hy\}$. It is easy to verify that we have defined a metric on $M/H$
which is compatible with the quotient topology. The inequality $w(M/H)\le w(M)$ follows from the fact that the quotient mapping $M\to M/H$
is open.
\epf

\blem
\label{l:maps}
Let $f:A\to M$ be a mapping from an infinite set $A$ to a metric space $(M,d)$. Suppose that either $w(M)<|A|$ or $M$ is compact.
Then for every $\delta>0$ there exist distinct elements $x,y\in A$ such that $d(f(x),f(y))<\delta$.
\elem

\bpf
The conclusion clearly is true if $f$ is not injective. Suppose that $f$ is injective, and let $\delta>0$ be given.
If  $w(M)<|A|$, the subset $f(A)$ of $M$ cannot be discrete, hence it contains a pair of distinct  $\delta$-close points.
Similarly, if $M$ is compact, the infinite set $f(A)$ cannot be closed discrete, so again there exists a pair of distinct
$\delta$-close points in $f(A)$, and the lemma follows.
\epf

\bprp\label{p:weight}
Let $G$ be a topological group.
Suppose that there exists an integer $n$ such that $w(K)<|\widehat G_n|$ for every compact subset $K$ of $G$.
Then  $1_G$ is not an isolated point in some $\widehat G_m\cup\{1_G\}$ with $1\leq m\leq n^2$.
\eprp

\begin{proof}
Observe that our assumption implies that $G$ and $\widehat G_n$ are infinite. Indeed,
for every finite subset $K$ of $G$ we have $|K|=w(K)<|\widehat G_n|$. If $|\widehat G_n|$ is finite,
it follows that so is $G$, and $|G|<|\widehat G_n|$. This is not possible, since a finite
group cannot have more irreducible unitary representations than its cardinality (\cite[Theorem 7]{Serre}).

Let $k=n^2$.
It suffices to prove that $1_G$ is not isolated in $\bigcup_{1\le m\le k} \widehat G_m\cup\{1_G\}$.
 Let $F$ be a compact subset of $G$ and $\gep>0$. We must prove that the neighborhood $W(1,F,\gep)$
of the trivial representation meets $\bigcup_{1\le m\le k} \widehat G_m$.

{\em Step 1.}
It suffices to construct  a non-trivial  irreducible unitary representation $\rho$ of $G$ of dimension $\le k$
which has an $(F,\gep)$-invariant unit vector $v$. Indeed, the function $f$ on $G$ defined by
$f(g)=(\rho(g)v,v)$ is a function of positive type associated with $\rho$, and for every $g\in F$ we have
$|f(g)-1|=|(\rho(g)v-v,v)|\le\norm{\rho(g)v-v}<\gep$. Thus
$[\rho]\in W(1,F,\gep)\cap \bigcup_{1\le m\le k} \widehat G_m\cup\{1_G\}$.


{\em Step 2.}
It suffices to construct  a $k$-dimensional unitary representation $\tau$ of $G$ without non-zero invariant vectors
which has an $(F,\gep/n)$-invariant unit vector $v$. Indeed, writing the representation as a direct sum of $s$ irreducible unitary representations
$\tau=\bigoplus\limits_{i=1}^s \tau_i$,
we can write $v=v_1+\dots+v_s$, where $s\le k$ and the vectors $v_i$ are $\tau_i$-invariant
and pairwise orthogonal. If $g\in F$, then $\tau(g)v-v=\sum_{i=1}^s (\tau_i(g)v_i-v_i)$ is a sum of orthogonal vectors,
hence $\norm{\tau_i(g)v_i-v_i}\le \norm{\tau(g)v-v}<\gep/n$ for every $i$.
Since $\sum_{i=1}^s \|v_i\|^2=1$,
one of the vectors $v_i$ has length $\ge1/\sqrt{s}\ge 1/\sqrt{k}=1/n$. Suppose it is the
vector $v_1$. Then $u=v_1/\|v_1\|$ is a unit vector which is $(F, \gep)$-invariant.
Indeed, for every $g\in F$ we have $\norm{\tau_1(g)u-u}=\norm{\tau_1(g)v_1-v_1}/\|v_1\|<\gep/(n\|v_1\|)\le\gep$.
We are in the situation of Step~1: $\tau_1$ is an irreducible representation of dimension $\le k$ for which $u$ is an
 $(F,\gep)$-invariant unit vector. Note that $\tau_1$ is non-trivial because $\tau$ has no non-zero invariant vectors.

We are going to construct a unitary representation $\tau$ with the properties described in Step 2.
The idea of the construction is the following.
Let $\mathcal F_n$ be the set of classes of all $n$-dimensional unitary representations of $G$ (which may be reducible).
Since $\widehat G_n\subseteq \mathcal F_n$ is ``big", we can find two distinct classes $[\rho_1], [\rho_2]\in \widehat G_n$
which are ``almost the same on $F$". We use $\rho_1$ and $\rho_2$ to construct a $k$-dimensional unitary representation $\tau$
on the space of all $n\times n$ complex matrices. We define $\tau$ by $\tau(g)A= \rho_1(g)A\rho_2\obr(g)$. Shur's Lemma
implies that $\tau$ has no non-zero invariant vectors, and the identity matrix $I$ will be $(F,\gep)$-invariant. We now elaborate.

{\em Step 3}.
Equip $\mathbb U(n)$ with any compatible bi-invariant metric $b$,
and equip $C(F,\, \mathbb U(n))$ with the sup-metric.
The compact group $ \mathbb U(n)$ isometrically acts on $C(F,\, \mathbb U(n))$ by conjugation. Denote by $O_n$ the orbit
space. According to Lemma~\ref{l:action}, the space $O_n$ has a natural metric compatible with the quotient topology.

Consider the set $\mathcal C$ of all continuous homomorphisms $G\to \mathbb U(n)$
(we do not consider any topology on $\mathcal C$). The group $\mathbb U(n)$ acts on $\mathcal C$ by conjugation, and the orbit space
is the set $\mathcal F_n$.
There is a natural map $\mathcal F_n\to O_n$ that is obtained from the restriction map $\mathcal C\to C(F,\, \mathbb U(n))$ by passing
to quotients. (We do not claim that the mapping  $\mathcal F_n\to O_n$ is continuous if $\mathcal F_n$ is equipped with the Fell topology.)
The image of $[\rho]\in \mathcal F_n$ in $O_n$ under this mapping will be denoted by $[\rho]|_F$. We have the following commutative diagram:
$$
\begin{CD}
\mathcal C @>>> C(F,\, \mathbb U(n))\\
@VVV@VVV\\
\mathcal F_n@>>>O_n
\end{CD}
$$
The horizontal arrows are restriction mappings, the vertical arrows are quotients by the action of $\mathbb U(n)$.

{\em Step 4}.
We claim that
for every $\delta>0$ we can find two homomorphisms $\rho_1, \rho_2:G\to \mathbb U(n)$
which are $\delta$-close
on $F$ and determine non-equivalent irreducible representations of $G$ on $\C^n$ (see \cite{ferher:dtopologies}).

Consider the mapping $[\rho]\to[\rho]|_F$ from $\widehat G_n\subseteq \mathcal F_n$ to $O_n$.
If $F$ is infinite,
we have $w(C(F,\, \mathbb U(n)))=w(F)<|\widehat G_n|$  \cite[Theorem 3.4.16]{Eng}
and $w(O_n)\le w(C(F,\, \mathbb U(n)))$ (Lemma~\ref{l:action}), so $w(O_n)<|\widehat G_n|$.
If $F$ is finite, $C(F,\, \mathbb U(n))$ and $O_n$ are compact. In either case we can apply Lemma~\ref{l:maps},
%
%
which implies that we can find distinct elements $[\rho_1], [\rho_2]\in \widehat G_n$ such that
$[\rho_1]|_F, [\rho_2]|_F$ are $\delta$-close in $O_n$. The distance between $[\rho_1]|_F$ and $[\rho_2]|_F$
in $O_n$ is equal to  
$$
\inf\{\hbox{dist\,}(\rho_1{}_{|F}, A\rho_2{}_{|F}A\obr): A\in \mathbb U(n)\},
$$
where dist refers to the distance in $C(F,\, \mathbb U(n))$. It follows that
replacing if necessary $\rho_2$ by an equivalent representation, we may assume that $\rho_1$ and $\rho_2$ 
are $\delta$-close on $F$,
as claimed.
%

{\em Step 5}.
Let $E=\End \C^n$ be the $n^2$-dimensional
Hilbert space of endomorphisms of $\C^n$. The scalar product on $E$ is given by the formula $(A,B)=\Tr(AB^*)/n$.
If $\rho_1,\rho_2$ are two irreducible unitary representations of $G$ on $\C^n$,
the formula $\tau(g)A= \rho_1(g)A\rho_2\obr(g)$ ($g\in G, \ A\in E$) defines a unitary representation
of $G$ on $E$ that does not contain non-zero invariant vectors. Indeed, if $A\in E$ is such that $\tau(g)A=A$ for all  $g\in G$,
then $A$ intertwines
$\rho_1$ and $\rho_2$ and hence is zero by Shur's Lemma.

If $I\in E$ is the identity mapping on $\C^n$, then $I$ is unit vector in $E$, since $(I,I)=\Tr(II^*)/n=1$.
By the continuity 
of the map $A\longmapsto  \Tr(A-I)(A-I)^*$,
there exists
$\delta>0$ such that the following holds: if $A\in \mathbb U(n)$ and $A$ is $\delta$-close to $I$ (that is, $b(A,I)<\delta$),
then
$$
\Tr(A-I)(A-I)^*<\gep^2/n.
$$
 Pick two homomorphisms $\rho_1, \rho_2:G\to \mathbb U(n)$ which are $\delta$-close
on $F$ and determine non-equivalent irreducible representations of $G$ on $\C^n$ (Step 4). Then $I\in E$ is an $(F,\gep/n)$-invariant
vector for the representation $\tau$ on $E$ constructed from $\rho_1, \rho_2$ as above. Indeed, let $g\in F$. Put
$A=\tau(g)I=\rho_1(g)\rho_2\obr(g)$.
Since the operators $\rho_1(g)$ and $\rho_2(g)$ are $\delta$-close and the metric $b$ is bi-invariant, $A$ is $\delta$-close to $I$.
Therefore we have
$$
\|\tau(g)I-I\|^2=\|A-I\|^2=\Tr(A-I)(A-I)^*/n<\gep^2/n^2
$$
and $\|\tau(g)I-I\|<\gep/n$. Thus $I$ is $(F,\gep/n)$-invariant.

We have constructed  a $k$-dimensional unitary representation $\tau$ of $G$ without non-zero invariant vectors
which has an $(F,\gep/n)$-invariant unit vector. According to Step~2, this completes the proof.

\end{proof}

The dual $\widehat G$ of a discrete group $G$ in general contains  infinite-dimensional irreducible unitary representations
(see \cite{Robert}). We denote by $\bigcup\limits_{n\in \N} \widehat G_n$ the subset of $\widehat G$ consisting
of all finite dimensional unitary representations. Applying Proposition \ref{p:weight} to the finite
subsets of a discrete group $G$, we obtain the following Corollary.
\bcor\label{co:prewang}
Let $G$ be a discrete group such that $1_G$ is an isolated point in $\bigcup\limits_{n\in \N} \widehat G_n$.
Then $\widehat G_n$ is finite for all $n\in \N$.
\ecor

Corollary \ref{co:prewang} yields the following result of Wang \cite{wang:1975}.
(Remark that property (T) is treated in the next section).
\bcor\label{co:wang}
If a discrete group $G$ has property (T), then $\widehat G_n$ is finite for all $n\in \N$.
\ecor
\bpf
If $G$ has property (T), then $1_G$ is isolated in $\widehat G\cup\{1_G\}$.
\epf

The {\em Bohr topology} on a topological group $G$ is the finest precompact group topology
that is coarser than the given topology on $G$. The {\em Bohr compactification} of $G$ is the completion
of $G$ equipped with the Bohr topology. Proposition~\ref{p:weight} also implies the following.

\bthm
\label{t:strong}
Let $G$ be topological group, $\kappa$ a cardinal such that $w(K)\le \kappa$
for every compact subset of $G$.
If  $1_G$ is an isolated point in $\widehat G_n\cup\{1_G\}$ for every $n$, then
$w(bG)\le\kappa$, where $bG$ is the Bohr compactification of $G$.
\ethm

\bpf
According to Proposition~\ref{p:weight}, $|\widehat G_n|\le\kappa$ 
for every $n$. The compact group $bG$ therefore
has $\le \kappa$ non-equivalent
irreducible unitary representations and hence has weight $\le \kappa$
by the Peter--Weyl theorem.
\epf

The case $\kappa=\omega$ of the previous theorem deserves to be explicitly stated:

\bcor
\label{c:strong}
Let $G$ be topological group such that every compact subset of $G$ is metrizable.
If  $1_G$ is an isolated point in $\widehat G_n\cup\{1_G\}$ for every $n$, then
the Bohr compactification of $G$ is metrizable.
\ecor

Note that this Corollary implies the following assertion which is stronger than Theorem~\ref{t:2}:
{\em if $G$ is a non-metrizable precompact group such that all compact subsets of $G$ are metrizable,
then $1_G$ is not an isolated point in $\widehat G$}.

Our proof of Theorem~\ref{t:3} is based on the following lemma. A  proof can be found e.g. in \cite{comfort:handbook}.
For the reader's convenience we provide the proof.

\blem
\label{l:omega}
Every compact group of weight $\omega_1$ has a dense countable subgroup.
\elem
\bpf
\A compact space $X$ is {\em dyadic} if there exists a mapping $2^\kappa\to X$
of a Cantor cube $2^\kappa$ onto $X$. All compact groups are dyadic (moreover, any compact $G_\delta$-subset of any
topological group is dyadic, see \cite{UspPrague, UspNotes}),
and any dyadic compact space of weight $\le\mathfrak c=2^\omega$ is separable, being an image of the separable space $2^{\mathfrak c}$.
\epf

We now prove the following stronger version of Theorem~\ref{t:3}:
{\em If $K$ is a non-metrizable compact group, then $K$ has a dense subgroup $G$ such that $1_G$ is not isolated
in $\widehat G$}.

Embedding $K$ into the product of unitary groups, we see that there exists a continuous homomorphism $f:K\to K'$ onto a compact group
of weight $\omega_1$. Let $G'$ be a dense countable subgroup of $K'$
(Lemma~\ref{l:omega}). Put $G=f\obr(G')$. Since $f$ is open, $G$ is dense in $K$. Consider the continuous mapping
$f^*:\widehat{G'}\to \widehat{G}$ dual to $f$.
Since $f$ is onto, $f^*$ is injective. Therefore,  $f^*$
sends $1_{G'}$ to $1_G$ and $\widehat{G'}\stm \{1_{G'}\}$ to $\widehat{G}\stm \{1_{G}\}$.
According to Theorem~\ref{t:2}, $1_{G'}$ is in the closure of $\widehat{G'}\stm \{1_{G'}\}$. It follows that
$1_{G}$ is in the closure of $\widehat{G}\stm \{1_{G}\}$. Equivalently, $1_{G}$ is not isolated in $\widehat{G}$.

\section{The property (T)}

The Kazhdan property (T) for topological groups was introduced in \cite{kazhdan}. This property
has several consequences on the structure of the groups that satisfy it. For $\sigma$-compact locally compact
groups Property (T)
is equivalent to the fixed point property for isometric affine actions on real Hilbert spaces
\cite[Theorem 2.12.4]{bhv:book}.
See \cite{bhv:book, dixmier, harp_vale} for further details on this important topic.

In the previous sections, we saw that for every metrizable precompact group $G$ the dual
$\widehat G$ is discrete (Theorem~\ref{t:1}). In contrast, we have the following result.

\bthm
\label{t:4}
If $G$ is an Abelian, countable precompact group, then $G$ does not have property (T).
\ethm

The result is no longer true if ``Abelian" is dropped. Indeed, certain compact Lie groups
admit dense countable subgroups which have property (T) as discrete groups
\cite[Theorem 6.4.4]{bhv:book} and hence
also as precompact topological groups.

We begin with the following characterization of property (T) for precompact groups.

\bthm\label{th:T}
Let $G$ be a dense subgroup of a compact group $H$. Then $G$ has property
(T) if and only if the following condition holds: there exist a compact subset $K\subseteq G$
and $\gep > 0$ such for every
function $f$ of positive type with zero integral over $H$ (with
respect to the Haar measure) the distance from $f|_{K}$ to $1$ in the space $C(K)$ (with
respect to the sup-norm) is $\geq \gep$.
\ethm
\bpf
Functions of positive type on $H$ are the functions of the form
$g\mapsto (\rho(g)v, v)$ associated with unitary representations $\rho$ of $H$ (here $v$ is a vector in the
Hilbert space $\mathcal H$ of the representation).
If a representation $\rho$ does
not contain the trivial representation (that is, there are no non-zero invariant vectors),
the integral of $(\rho(g)v, v)$ over $H$ is zero, for every $v\in \mathcal H$ (this follows, for example,
from orthogonality relations). Conversely, if the integral of
$f(g) = (\rho(g)v, v)$ is zero, then $v$ lies in $\mathcal H_1 = \mathcal H_2^\perp$,
where $\mathcal H_2\subseteq \mathcal H$ is the subspace of invariant vectors.
Therefore $f$ is associated with a representation on $\mathcal H_1$ without non-zero invariant vectors.

Now suppose that $G$ has property (T). Let $(K,\delta)$ be a Kazhdan pair for $G$
and assume, without loss of generality, that $e\in K$ .  Suppose that $f$ is a function
of positive type on $H$ with zero integral over $H$. Then $f(g) = (\rho(g)v, v)$ for some unitary representation $\rho$ of
$H$ on a Hilbert space $\mathcal H$ without non-zero invariant vectors  and some $v\in \mathcal H$. We first consider the case
when $v$ is a unit vector, that is, $f(e)=1$, where $e$ the neutral element of $G$.
By the definition of a Kazhdan pair, there exists $g\in K$ such that $\norm{\rho(g)v-v}\ge \delta$. Writing $f(g)=(\rho(g)v, v)=a+bi$,
we have $\delta^2\le \norm{\rho(g)v-v}^2=2-2a$ and thus $a\le 1-\delta^2/2$.
It follows that $|f(g)-1|\ge1-a\ge \delta^2/2$.

In the general case (when we no longer assume that $f(e)=1$), we replace $f$ by $f/f(e)$ (we exclude the trivial case $f(e)=0$).
If there is no $\gep$ such that the pair $(K,\gep)$ satisfies the conditions of the theorem, there is a sequence $(f_n)$ of functions
of positive type with zero integral over $H$ such that the sequence $(f_n|_K)$ uniformly converges to 1 on $K$. We may assume,
without loss of generality, that $f_n(e)=1$ for every $n$ (otherwise replace $f_n$ by $f_n/f_n(e)$). According to the last inequality
of the previous paragraph, the distance from $f_n|_K$ to 1 in $C(K)$ is at least $\delta^2/2$. This contradicts the fact that
the sequence $(f_n|_K)$ uniformly converges to 1 on $K$ and proves the existence of a number $\gep>0$ such that the pair
$(K,\gep)$ satisfies the conditions of the theorem.

Conversely, suppose that $G$ does not have property (T). We must prove that for every compact subset $K\subseteq G$
and $\gep>0$ there exists a function $f$ of positive type with zero integral over $H$ such that $|f(g)-1|< \gep$ for every $g\in K$.
Since $(K,\gep)$ is not a Kazhdan pair, there exist a unitary representation $\rho$ of $G$ on a Hilbert space $\mathcal H$
without invariant vectors and a unit vector $v\in \mathcal H$ which is $(K,\gep)$-invariant, that is, $\norm{\rho(g)v-v}<\gep$
for every $g\in K$. The representation $\rho$ of $G$ can be extended to a unitary representation $\overline{\rho}$ of $H$ on $\mathcal H$.
Since $\rho$ admits no invariant vectors, the same is true for $\overline{\rho}$. Let $f$ be the function on $H$ defined by
$f(g)=(\overline{\rho}(g)v, v)$. The function $f$ is as required: it is of positive type with zero integral over $H$, and for every $g\in K$
we have $|f(g)-1|=|(\rho(g)v-v,v)|\le \norm{\rho(g)v-v}\cdot\norm{v}<\gep$.
\epf

\bcor\label{co:T}
Let $G$ be a dense subgroup of a compact group $H$. If $G$ has property
(T), then there exist a compact subset $K\subseteq G$, a real (signed) measure $\mu$ on $K$ and
a real number $c$ such that for every real function $f$ of positive type on $H$ with $f(e) = 1$ and zero integral
over $H$ we have $$\int_K f\,\mu < c < \int_K 1\,\mu.$$
\ecor
\bpf
Consider the convex set $P$ of all real functions $f$ of positive type on $H$ such that $\int_H f=0$ and $f(e) = 1$.
Let $K\subseteq G$ be a compact subset with the property described in Theorem~\ref{th:T}.
According to the theorem, the image of $P$ under the restriction map $f\mapsto f|_{K}$ is at a positive
distance from $1$. The Hanh-Banach theorem implies that $1$ can be separated from
the convex set $P|_{K}$ by a linear functional.
\epf

Following Lubotzky and Zimmer \cite{lub_zim:89}, if $\mathcal R$ is a set of classes of unitary representations of $G$,
we say that $G$ has \emph{Property (T) with respect to $\mathcal R$} 
if $1_G$ is isolated in $\mathcal R\cup \{1_G\}$.

\brem\label{rm:T}
Let $G$ be a precompact group,  and let $G_d$ denote the same group with the discrete topology.
Let $\mathcal R$ be the set
of equivalence classes of finite-dimensional unitary representations of $G$ without non-zero invariant vectors.
It is readily seen that the proof of Theorem~\ref{th:T} shows that $G_d$ has property (T) with respect to $\mathcal R$
if and only if there exist a finite subset $K\subseteq G$ and $\gep >0$ such that for every
function $f$ of positive type with zero integral over the completion of $G$ the distance from $f|_K$ to $1$
in the space $C(K)$ (with respect to the sup-norm) is $\geq \gep$.
\erem

\bthm\label{t:trick}
Let $G$ be a precompact group such that every compact subset of $G$ is countable, and let
$G_d$ be the same 
group equipped with the discrete topology.
Let $\mathcal R$ be the set
of equivalence classes of finite-dimensional unitary representations of $G$ without non-zero invariant vectors.
The following assertions are equivalent:
\begin{enumerate}
\item[(a)] $G$ has property (T);

\item[(b)] $G_d$  has property (T) with respect to $\mathcal R$.
\end{enumerate}
\ethm
\bpf
(b) $\Rightarrow$ (a) is clear.

(a) $\Rightarrow$ (b).
Let $H$ be the completion of $G$. Suppose that $G$ has property (T) but $G_d$ does not have property (T)
with respect to $\mathcal R$.
Let $P$ be the convex set that was used in the proof of Corollary~\ref{co:T}: $P$ is the set of all real-valued
functions $f$ of positive type on $H$ with zero integral and such that $f(e)=1$.
Using Corollary~\ref{co:T}, find a compact set $K\subseteq G$ such that $1|_{K}$ is not in the weak closure
of the convex set $P|_{K}$ in the Banach space $C(K)$ (the weak closure and the norm closure are the same
because of convexity, see \cite[Theorem 3.12]{rudin}). Since all compact subsets of $G$ are countable, we
can write $K$ as the set of points of a sequence $(x_n)$.
According to Remark \ref{rm:T},
for every finite subset $F$ of $G$ and $\gep >0$ there exists $f\in P$
such that the distance from $f|_{F}$ to $1|_{F}$
in the space $C(F)$ is $< \gep$. For every positive integer $n$ apply this to $F=\{x_i:i\le n\}$
and $\gep=1/n$.
We obtain a function $f_n\in P$ such that $|f_n(x_i)-1|<1/n$ for every $i\le n$.
The sequence $(f_n)$ pointwise converges
to $1$ on $K$. All functions in $P$ are uniformly bounded by 1,
so Lebesgue's dominated
convergence theorem implies that the sequence $(f_n|_K)$ weakly converges to $1$ in $C(K)$.
Therefore, $1|_K$ is in the weak closure of $P|_{K}$, contrary to our assumption.
\epf

%
%

We are now in a position to prove Theorem~\ref{t:4}:

{\em If $G$ is an Abelian, countable precompact group, then $G$ does not have property (T).}

\bpf
Let $H$ be the completion of $G$. By Theorem~\ref{t:trick}, we must verify that $G_d$
does not have property (T) with respect to $\mathcal R$, where $\mathcal R$ is the same as in Theorem~\ref{t:trick}.
It suffices to prove that for every finite set $F\subseteq G$ and $\gep>0$ there is a continuous character $\chi:H\to \mathbb U(1)$,
$\chi\ne1$, such that $|\chi(x)-1|<\gep$ for every $x\in F$.

Assume the contrary. Let $F\subseteq G$ be a finite set such that for some $\gep>0$ and every $\chi\in \widehat H\stm\{1\}$
the restriction $\chi|_{F}$ is at the distance $\ge\gep$ from $1|_{F}$. Consider the restriction homomorphism
$\phi: \widehat H\to \mathbb U(1)^F$ defined by $\phi(\chi)=\chi|_{F}$.
According to our assumption, $\phi$ has a trivial kernel and hence is injective.
Also, $1|_{F}$ is at the distance $\ge\gep$ from $\phi(\widehat H\stm\{1\})$. It follows that $1|_{F}$ is an isolated point
in the group $\phi(\widehat H)$. Since every topological group which has an isolated point is discrete, we conclude that
$\phi(\widehat H)$ is a discrete subgroup of $\mathbb U(1)^F$. This is a contradiction, since the compact group
$\mathbb U(1)^F$ cannot have infinite discrete subgroups.
\epf

\section{Questions}

\bque\label{q:1}
Does there exist a non-compact precompact Abelian group with property (T)?
\eque

By Theorem~\ref{t:4} such a group must be uncountable.

The definition of the Fell topology, given in \cite{bhv:book} and used in the present paper, is not the same
as the definition given in \cite{Marg}. Using the notation of Section~\ref{s:prelim}, the difference is the following.
In \cite{Marg}, functions $f_1, \dots, f_n\in P_\rho'$ of positive type are approximated by
functions $g_1,\dots, g_n\in P_\s'$ of positive type rather than by their sums in $P_\rho$ that we allowed. For locally compact groups $G$
the two definitions of the Fell topology on $\widehat G$ agree, according to \cite[Proposition F.1.4]{bhv:book}.

\bque
Do the two definitions of the Fell topology coincide on $\widehat G$: (a) for every topological group?
(b) for every precompact group?
\eque

\section{Acknowledgement}

We thank the referee for many valuable suggestions.

\end{document}